\input amstex
\documentstyle{amsppt}  
\mag=\magstep1  
\NoRunningHeads  
\NoBlackBoxes 
\define\id{\operatorname{id}} 
\define\diam{\operatorname{diam}} 
\define\ps{\operatorname{ps}} 
\define\dist{\operatorname{dist}} 
\topmatter  
\title Pseudo-orbits, pseudoleaves and geometric entropy of foliations\endtitle  
\author Andrzej Bi\'s and Pawe\l \/ Walczak\endauthor  
\address Instytut Matematyki, Uniwersytet \L\'odzki, Banacha 22, 
90238 \L\'od\'z, Poland\endaddress  
\email andbis\@ krysia.uni.lodz.pl and pawelwal\@ krysia.uni.lodz.pl \endemail  
\date November 4, 1996\enddate
\thanks The second author was supported by KBN grant 1037/P03/96/10.\endthanks 
\abstract We show that the entropy of a finitely generated pseudogroup (resp.,
of a foliation on a compact Riemannian manifold) can be calculated by suitable
counting separated pseudo-orbits (resp., pseudoleaves).\endabstract
\subjclass 53 C 12, 57 R 30, 54 H 20\endsubjclass 
\endtopmatter  
 
\document 
\subhead 0. Introduction\endsubhead 
Since several years, pseudo-orbits play an important role in the theory 
of classical  
dynamical systems. In particular, it is shown ([{\bf Mi}] and [{\bf BS}]) that
pseudo-orbits can be used to calculate the topological entropy 
of transformations. 
More recently, the similar result was obtained [{\bf Hur}] for the 
inverse-image 
entropy introduced earlier by Remi Langevin and the second author [{\bf LW1}]. 
 
On the other hand, geometric entropy $h(\Cal F)$ of a foliation $\Cal F$ 
of a compact  
Riemannian manifold has been introduced [{\bf GLW}] and shown to be a handful  
tool to study topology and dynamics of foliated manifolds 
( [{\bf Hu1 - 3}], [{\bf GLW}], [{\bf  LW2}], [{\bf Eg}], [{\bf IT}], etc.). 
The entropy $h(\Cal F)$ can be  
calculated either by counting the number of points separated along the 
leaves or by counting the number of separated orbits of holonomy pseudogroups 
generated by nice  
coverings by charts distinguished by $\Cal F$. Also, pseudoleaves of foliations have been
defined by Takashi Inaba [{\bf In}] who has shown (among the other results) 
that expansive C$^1$-foliations (in the sense of [{\bf IT}]) of codimension-one
which have the pseudoleaf tracing property are topologically stable. 
 
In this article, we show that the geometric entropy $h(\Cal F)$ of any 
foliation $\Cal F$ of a compact Riemannian manifold $M$ coincides with that 
calculated by suitable counting separated pseudoleaves (Theorem 2 in Section 
4). To this end, we study the entropy of finitely generated pseudogroups of 
local transformations of compact metric spaces and  
we show that it can be calculated by counting (again, in a suitable way) 
separated pseudo-orbits (Theorem 1 in Section 2). In Section 5, we sketch an 
easy proof of Theorem 2 for foliated bundles. Finally, in Section 6 we provide
an example of a group acting on $S^1$ for which the "usual" formula for the 
entropy in terms of pseudo-orbits does not work: It gives the number strictly
bigger than the entropy calculated in terms of separated orbits.

Some significance of these results could be observed if one would try to 
calculate (or, to estimate)  entropies of  pseudogroups or foliations 
with the aid of computers.  

The authors are grateful to Takashi Inaba for his valuable comments made during
his visit to the University of \L\'od\'z. 

\subhead 1. Pseudo-orbits of pseudogroups\endsubhead 
Let $(X, d)$ be a compact metric space and $G$ be a finitely generated 
pseudogroup  of local transformations of $X$. For any $g\in G$, $U_g\subset X$ 
will denote the domain of $g$. Let us fix a finite symmetric generating set  
$G_1\subset G$. Recall that "symmetric" means that $e = \id_X\in G_1$ 
and $G_1^{-1}\subset G_1$. For any $k\in\Bbb N$ let $G_k = (G_1)^k = \{ 
(h_1,\dots h_k); h_i\in  G_1\}$. We shall identify $G_k$ with a subset of 
$G_{k+1}$ consisting of all the sequences of the form $(e, h_1,\dots h_k)$, 
$h_i\in G_1$. Also, given two sequences $h = (h_1,\dots h_k)\in G_k$ and 
$g = (g_1\dots g_l)\in G_l$ we shall denote by $(h, g)$ the element 
$(h_1,\dots h_k, g_1,\dots g_l)\in G_{k+l}$. Finally, for any $h =  
(h_1,\dots k_k)\in G_k$ we shall denote by the same symbol $h$ the map  
$h_1\circ\dots\circ h_k\in G$. Therefore, we may write that $\{ e\} = 
G_0\subset G_1\subset G_2\subset\dots\subset G_\infty = \cup_kG_k$.  
 
The entropy of the system $(G, G_1)$ can be defined whenever the maps of 
$G$ are homeomorphisms between open sets. However, in general it can be 
infinite and the most natural setting ensuring its finiteness is the 
category of locally Lipschitz maps. Therefore, throughout the paper we 
shall assume that the maps generating pseudogroups are Lipschitz.  
 
Let us fix $\alpha\ge 0$. An $\alpha$-{\it pseudo-orbit} of $G$ 
is a map $x:D_x\to X$ such that $D_x\subset G_\infty$, $e\in D_x$, 
for any $h\in G_1$ and $g\in D_x$, $(h, g)\in D_x$ iff $x(g)\in U_h$, 
$U_h$ being the domain of $h$, and then 
$$d(h(x(g)), x(h, g))\le\alpha .\tag1$$ 
 
If $G$ is a group of global transformations, then $D_x = G_\infty$ for any  
pseudo-orbit $x$ and inequality (1) should be satisfied for all $g\in G_\infty$ 
and $h\in G_1$. In particular, if $G = \Bbb Z$ acts on $X$ and $f$ is 
a generator, then $x:\Bbb Z\to X$ is just a sequence of points and (1) 
reduces to 
$$d(f(x_n), x_{n+1})\le\alpha ,\quad n\in\Bbb Z$$ 
i. e. $x$ becomes a pseudo-orbit of $f$ in the classical sense [{\bf Bo}]. 
 
Let us denote by $Y_\alpha$ the space of all the $\alpha$-pseudo-orbits. 
Clearly, $Y_\alpha\subset Y_\beta$ whenever $\alpha\le\beta$ and $Y_0 = 
\cap_{\alpha\ge 0}Y_\alpha$ is the space of "real" $G$-orbits. Recall, 
that $x\in Y_0$ iff $x = x_p$ for some $p\in X$, $D_{x_p} = \{ 
g\in G_\infty ; p\in U_g\}$ and $x_p(g) = g(p)$ for all $g\in D_{x_p}$.  
 
For any pseudo-orbits $x$ and $y$ in $Y_\alpha$ put 
$$d_0(x, y) = \sum_{k=0}^\infty\frac{1}{2^k}\max\{ d(x(g), y(g)); 
g\in G_k\cap D_x\cap D_y\}\tag2$$ 
and 
$$d_1(x, y) = \inf\{\sum_{i=0}^{m-1} d_0(z_i, z_{i+1}); z_i\in Y_\alpha , 
z_0 = x, z_m = y\ \text{and}\ m\in\Bbb N\} .\tag3$$ 
 
\proclaim{Lemma 1} For any $\alpha$, $(Y_\alpha , d_1)$ is a metric space. 
\endproclaim 
 
\demo{Proof} The symmetry and triangle inequality for $d_1$ as well as the 
relation $d_1(x, x) = 0$ for all $x\in Y_\alpha$ follow immediately from (2) 
and (3). It remains to show that the condition $d_1(x, y) = 0$ implies 
the equality $x = y$ which is equivalent to the following: 
$$D_x\cap G_k = D_y\cap G_k\ \text{and}\ x = y\ \text{on}\ D_x\cap D_y\cap
G_k\tag4$$ 
for $k = 0, 1, 2, \dots $. 
 
First, observe that 
$$\sum_{i=0}^{m-1} d_0(z_i, z_{i+1})\ge\sum_{i=0}^{m-1} d(z_i(e), 
z_{i+1}(e))\ge d(x(e), y(e))$$ 
whenever $z_0 = x$ and $z_m = y$. Therefore, the condition $d_1(x, y) = 0$ 
implies that $x(e) = y(e)$ and then (4) is satisfied with $k = 0$.  
 
Next, observe that the equality $x(e) = y(e) = p$ implies that $D_x\cap G_1 
= D_y\cap G_1 = G(p)\cap G_1$, where $G(p) = \{ h\in G; p\in U_h\}$. Find 
$\epsilon > 0$ such that $B(p, \epsilon )\subset U_h$ for all $h\in G(p)\cap 
G_1$. If $z_0, z_1, \dots , z_m\in Y_\alpha$, $z_0 = x$, $z_m = y$ and 
$\sum_i d_0(z_i, z_{i+1}) < \frac{\epsilon}{2}$, then 
$$d(p, z_j(e))\le\sum_{i\le j} d(z_i(e), z_{i+1}(e))\le \sum_{i\le j} 
d_0(z_i, z_{i+1}) < \frac{\epsilon}{2}$$ 
and $z_j(e)\in U_h$ for any $h\in G(p)\cap G_1$ and $j\le m$. Therefore, 
$$d(x(h), y(h))\le\sum_i d(z_i(h), z_{i+1}(h))\le 2\cdot\sum_i d_0(z_i, 
z_{i+1}) < \epsilon$$ 
for all $h\in G(p)\cap G_1$. This shows that the condition $d_1(x, y) = 0$ 
implies (4) with $k = 1$. 
 
Proceed by induction to complete the proof. \qed 
\enddemo 
 
\proclaim{Lemma 2} The space $(Y_\alpha , d_1)$ is compact. 
\endproclaim 
 
\demo{Proof} Take any sequence $(x_n)$ of $\alpha$-pseudo-orbits. Choose a  
subsequence $(x^0_n)$ of $(x_n)$ such that $x^0_n(e)\to p\in X$. Put 
$x(e) = p$ and $D_x\cap G_1 = G(p)\cap G_1$. Find a subsequence $(x^1_n)$ 
of $(x^0_n)$ such that $x^1_n(h)\to p_h\in X$ for any $h\in D_x\cap G_1$. 
Put $x(h) = p_h$ and $D_x\cap G_2 = \{ (h', h); h\in D_x\cap G_1, 
h'\in G_1\ \text{and}\ p_h\in U_{h'}\}$. Proceed by induction to define 
$D_x\cap G_k$, to find subsequences $(x^k_n)$ such that $x^k_n(g)\to 
p_g\in X$ for all $g\in D_x\cap G_k$, $k\in\Bbb N$,  
and put $x(g) = p_g$ for all such $g$. 
 
If $h\in D_x\cap G_1$, then $h\in D_{x^0_n}$ and  $d(x^0_n(h),  
h(x^0_n(e))\le\alpha$ for all $n$ large enough. Therefore, 
$d(x(h), h(x(e)))\le\alpha$.  
Also, if $(h', h)\in D_x\cap G_2$, then $h'\in D_{x^1_n}$ and  $d(x^1_n(h', h),
h'(x^1_n(h)))\le\alpha$ for all $n$ large enough. Therefore, $d(x(h', h),  
h'(x(h)))\le\alpha$. Proceed by induction to show that $x\in Y_\alpha$. 
 
Obviously, $x^n_n\to x$. \qed 
\enddemo 
 
Now, take $g_0\in G_\infty$ and $x\in Y_\alpha$ such that $g_0\in D_x$. Define 
$y = \sigma_{g_0}(x)\in Y_\alpha$ by 
$$y(g) = x(g, g_0)$$ 
on $D_y = \{ g\in G_\infty ; (g, g_0)\in D_x\}$. It is easy to see that this 
definition is correct and that all the maps $\sigma_h$, $h\in G_1$, generate 
a pseudogroup, denoted by $G$ again, of local transformations of $Y_\alpha$.  
 
Note that 
$$\aligned d_0(\sigma_{g_0}x, \sigma_{g_0}x') &=  
\sum_k\frac{1}{2^k}\max\{d(x(g, g_0), x'(g, g_0)); g\in G_k\cap D_x\cap  
D_{x'}\}\\ &\le \sum_k\frac{1}{2^k}\max\{d(x(g), x'(g)); g\in G_{k+n}\cap 
D_x\cap D_{x'}\}\le 2^nd_0(x, x')\endaligned$$ 
and 
$$\aligned d_1(\sigma_{g_0}x, \sigma_{g_0}x') &= \inf\{\sum_{i\le m-1} 
d_0(z_i, z_{i+1}); z_0 = \sigma_{g_0}x\ \text{and}\ z_m = 
\sigma_{g_0}x')\}\\&\le\inf\{ \sum_{i\le m-1}d_0(\sigma_{g_0}z_i, 
\sigma_{g_0}z_{i+1}); z_0 = x\ \text{and}\ z_m = x'\}\le 2^nd_1(x, 
x')\endaligned$$ 
whenever $g_0\in G_n$. Therefore, $G$ acts on $Y_\alpha$ also {\it via} locally
Lipschitz transformations. 
 
\subhead 2. Entropy of pseudogroups\endsubhead 
As before, let $G$ be a pseudogroup on $X$ equipped with a finite symmetric  
generating set $G_1$. Following [{\bf GLW}] we shall say that two points $p$ 
and $q$ of $X$ are $(n, \epsilon)$-{\it separated} w. r. t. $G$
($n\in\Bbb N$, $\epsilon > 0$) if there exists  
$g\in G_n$ such that $\{ p, q\}\subset U_g$ and $d(g(p), g(q))\ge\epsilon$. 
A set $A\subset X$ is $(n, \epsilon)$-separated when any two points $p, 
q\in A$, $p\ne q$, have this property. Since $X$ is compact, separated 
sets are always finite and we may put 
$$N(n, \epsilon , X) = \max\{\# A; A\subset X\ \text{is $(n, 
\epsilon)$-separated}\},\tag5$$ 
$$N(\epsilon , X) = \limsup_{n\to\infty}\frac{1}{n}\log N(n, \epsilon , 
X)\tag6$$ 
and 
$$h(G, G_1, X) = \lim_{\epsilon\to 0^+} N(\epsilon , X).\tag7$$ 
The number $h(G, G_1, X)$ is the {\it entropy} of $G$ (w. r. t. $G_1$) on $X$.
 
\proclaim{Lemma 3} $h(G, G_1, X) = h(G, G_1, Y_0)$, $Y_0$ being the space of  
$G$-orbits. 
\endproclaim 
 
\demo{Proof} The inequality "$\le$" follows immediately from the estimate 
$$d_1(\sigma_g x, \sigma_g y)\ge d(x(g), y(g)) = d(g(p), g(q))$$ 
which holds for the orbits $x$ and $y$ of any points $p$ and $q$ of $X$ and any  
element $g\in G_\infty$ such that $\{ p, q\}\subset U_g$.
 
To prove "$\ge$" take two $(n, \epsilon )$-separated orbits $x$ and $y\in Y_0$ 
and find $g\in G_n$ such that $\{ x(e), y(e)\}\subset U_g$ and 
$d_1(\sigma_g(x), \sigma_g(y))\ge\epsilon$. Then also $d_0(\sigma_g(x), 
\sigma_g(y))\ge\epsilon$. Next,  
find $l(\epsilon )\in\Bbb N$ such that $\sum_{k > l(\epsilon )}2^{-k} <  
\epsilon/(2\diam X)$. Thus,  
$$d(g'g(x(e)), g'g(y(e))) = d(x(g',g), y(g',g)) = d(\sigma_gx(g'),  
\sigma_gy(g'))\ge\frac{\epsilon}{2}$$ 
for some $g'\in G_{l(\epsilon )}$ and the points $x(e)$, $y(e)$ are 
$(n + l(\epsilon ), \epsilon/2)$-separated on $X$. Consequently, 
$$N(n, \epsilon , Y_0)\le N(n + l(\epsilon ), \frac{\epsilon}{2}, X).$$ 

Pass to suitable limits to complete the proof. \qed
\enddemo

\proclaim{Lemma 4} For any $n\in\Bbb N$ and any $\epsilon > 0$ there exists 
$\alpha > 0$ such that $N(n, \frac{\epsilon}{3}, Y_0)\ge N(n, \epsilon  
, Y_\alpha )$.
\endproclaim 
\demo{Proof} Fix $n\in\Bbb N$ and $\epsilon > 0$. Let $L$ be the Lipschitz
constant for all the generators of $G$. Let $A = \{ x_i, i\le k\}$ be an
$(n, \epsilon )$-separated subset of $Y_\alpha$. For any $i$, let $x'_i$ be
the orbit of $x_i(e)$. An easy induction shows that
$$d(x_i(g), x'_i(g))\le\alpha\cdot (1 + L + \dots + L^{n-1})$$
for any $i$ and any $g\in G_n$. Therefore, for any $i$ and $j\le k$ there 
exists $g\in G_n$ such that
$$\align d(x'_i(g), x'_j(g))&\ge d(x_i(g), x_j(g)) - d(x_i(g), x'_i(g)) -
d(x_j(g), x'_j(g))\\&\ge\epsilon - 2\alpha\cdot (1 + L + \dots +
L^{n-1})\ge\frac{\epsilon}{3}\endalign$$
whenever $\alpha$ is small enough to satisfy $\alpha\cdot (1+L+\dots +L^{n-1})
\le\frac{\epsilon}{3}$. For such $\alpha$, $A' = \{ x'_i, i\le k\}\subset Y_0$
is $(n, \epsilon/3)$-separated.\qed
\enddemo

Now, let us say that $\alpha$-pseudo-orbits $x$ and $y$ are 
$(n, \epsilon )$-{\it strongly separated}
just when $d(x(g), y(g))\ge\epsilon$ for some $g\in G_n\cap D_x\cap D_y$. Let 
us denote by $N_\alpha (n, \epsilon )$ the maximal cardinality of a set 
$A\subset Y_\alpha$ such that any two elements $x, y\in A$, $x\ne y$,  are
$(n, \epsilon)$-separated in this sense. Note that in general 
$N_\alpha (n, \epsilon)\ne N(n, \epsilon , Y_\alpha)$ since strong separation 
of $\alpha$-pseudo-orbits in the above sense is not 
entirely equivalent to separation of points of $Y_\alpha$ by the action 
of $G$. For any sequence $\beta = (\beta_n)$ such that $0 < 
\beta_n\searrow 0$ we shall put 
$$N_\beta (\epsilon ) = \limsup_{n\to\infty}\frac{1}{n}\log N_{\beta_n} (n,
\epsilon),\tag8$$ 
$$N_{\ps} (\epsilon ) = \inf_{\beta}N_\beta (\epsilon )\tag9$$ 
and 
$$h_{\ps} (G, G_1, X) = \lim_{\epsilon\to 0^+}N_{\ps}(\epsilon ).\tag10$$ 
Our main result (Theorem 1) in this Section says that this "pseudo-entropy" 
$h_{\ps}$ coincides always with the "real " entropy $h$. 
 
\proclaim{Lemma 5} $N_\alpha (n, \epsilon )\le N(n, \epsilon , Y_\alpha)$
for all $n\in\Bbb N$, $\epsilon > 0$ and $\alpha\ge 0$.
\endproclaim 
 
\demo{Proof} If pseudo-orbits $x, y\in Y_\alpha$ are $(n, \epsilon )$-strongly 
separated and $d(x(g), y(g))\ge\epsilon$ for some $g\in G_n\cap D_x\cap D_y$, 
then 
$$d_1(\sigma_gx, \sigma_gy)\ge d(\sigma_gx(e), \sigma_gy(e)) = d(x(g),  
y(g))\ge\epsilon ,$$ 
i. e. $x$ and $y$ are separated under the action of $G$ on $Y_\alpha$.\qed 
\enddemo

\proclaim{Theorem 1} For any pseudogroup $G$ on $X$ and any finite symmetric  
generating set $G_1$ the equality 
$$h_{\ps} (G, G_1, X) = h(G, G_1, X)\tag11$$ 
holds. 
\endproclaim 
 
\demo{Proof} The inequality "$\ge$" in (11) is obvious since $Y_0\subset 
Y_\alpha$ for all $\alpha\ge 0$ and the notions of separation of points of 
$X$ under the action of $G$ and strong separation of orbits just coincide.
 
To prove "$\le$" observe that, by Lemmas 4 and 5, 
$$\align N_{\ps} (\epsilon )\le N_\beta (\epsilon)&\le 
\limsup_{n\to\infty}\frac{1}{n}\log N(n, \epsilon , Y_{\beta_n})\\ &\le
\limsup_{n\to\infty}\frac{1}{n}\log N(n, \frac{\epsilon}{3}, Y_0)
= N(\frac{\epsilon}{3}, Y_0)\endalign$$
for any sequence $\beta = (\beta_n)$ such that $\beta_n\le\epsilon\cdot 
(1+L+\dots +L^{n-1})^{-1}$, $n\in\Bbb N$ and $L$ being - as before - the 
Lipschitz constant for all the generators of $G$. Therefore, 
$$h_{\ps} (G, G_1, X)\le h(G, G_1, Y_0) = h(G, G_1, X)$$ 
by Lemma 3. \qed 
\enddemo

\demo{Remark} The proof above shows that
$$h(G, G_1, X) = \lim_{\epsilon\to 0^+}\lim_{n\to\infty}\frac{1}{n}
\log N_{\epsilon\alpha_n}(n, \epsilon)$$
with $\alpha_n = (1 + L + \dots + L^n)^{-1}$.
\enddemo

\subhead 3. Pseudoleaves and holonomy pseudo-orbits\endsubhead 
First, assume that $E$ is a vector space with an inner product $\langle\cdot ,
\cdot\rangle$. For any $q\ge 1$ and codimension-$q$ subspaces $E_1$ and $E_2$ 
of $E$ put 
$$\dist (E_1, E_2) = \max\{\min\{\| v_1 - v_2\| ; v_2\in E_2, \| v_2\| = 1\}; 
v_1\in E_1, \| v_1\| =1\} .\tag12$$ 
Clearly, $\dist$ is a distance function on the set of all codimension-$q$ 
subspaces of $E$ which is nearly proportional to the angle between subspaces.
 
Let $\Cal F$ be a C$^r$-foliation, $r\ge 1$, of a compact Riemannian manifold 
$(M, \langle\cdot , \cdot\rangle )$. Following [{\bf In}] we shall say that 
a complete submanifold $N$ of $M$ is an $\alpha$-pseudoleaf when $\dim N 
= \dim\Cal F$ and 
$$\dist (T_pN, T_p\Cal F)\le\alpha\tag13$$ 
for all $p\in N$.  
 
Next, take a nice cover $\Cal U$ of $(M, \Cal F)$ [{\bf HH}] and denote 
by $H_{\Cal U}$ the holonomy pseudogroup of $\Cal F$ acting on the space 
$X = X_{\Cal U}$ of plaques of the closed charts $\bar U$, $U\in\Cal U$. 
Since $\Cal U$ is finite, $X$ is compact.  
 
The space $X$ can be equipped with some natural distance functions $d = d_{\Cal
U}$ induced from the Riemannian structure of $M$. For example, one can equip 
$M$ with the Riemannian distance function $d_M$ and define $d_{\Cal U}$ as the  
Hausdorff distance between the closures of the plaques. Or, one can fix a 
complete transversal $T = T_{\Cal U}$ intersecting each plaque exactly once, 
identify $X$ with $T$ and equip $T$ with the Riemannian distance function 
coming from the induced Riemannian structure. All these choices are equivalent 
in the sense that any two metrics described above are quasi-isometric and 
therefore the entropies of $H_{\Cal U}$  
obtained from these metrics are equal ([{\bf GLW}], Prop. 2.6). So, let us 
fix such a "good" metric $d$ on $X$ once for ever and use it hereafter all the
time. 
 
Also, we consider the "standard" generating set $(H_{\Cal U})_1$ for 
$H_{\Cal U}$. It consists of all the elementary holonomy maps $h_{U, V}$ 
corresponding to overlapping charts $U, V\in\Cal U$ and defined by the 
following: 
$$h_{U, V}(P) = Q\ \text{iff the plaques $P\subset U$ and $Q\subset V$  
overlap.}\tag14$$ 
Since this is the only generating set to be considered, we shall omit the 
symbol $(H_{\Cal U})_1$ in all the formulae below. 
 
Now, observe that the space of pseudoleaves embeds naturally into the space of  
pseudo-orbits of  some holonomy pseudogroups $H_{\Cal U}$. To produce such an  
embedding fix a complete transversal $T_{\Cal U}$ and a subbundle $E$ of $TM$  
complementary to $T\Cal F$ and such that $E_p = T_pT_{\Cal U}$ for any $p\in  
T_{\Cal U}$. Also, for any two overlapping charts $U, V\in\Cal U$ fix another  
distinguished chart $W = W_{U, V}$ containing $\bar U\cup\bar V$ and assume 
that $\Cal U$ is tiny enough for all the plaques of the charts $W_{U, V}$ to be  
strongly convex subsets of the leaves. 
 
\proclaim{Lemma 6} There exist positive constants $C, C'$ and $\alpha_0$ such 
that for any $\alpha\le\alpha_0$, any $\alpha$-pseudoleaf $N$ and any $p\in N$ 
there exists an unique $(C\alpha)$-pseudo-orbit $x$ of $H_{\Cal U}$ such that 
$x(e) = p$, $x(D_x)\subset N$ and $d_N(x(h,g), x(g))\le C'd_{\Cal F}(h(x(g)), 
x(g))$ for any $g\in D_x$ and any generator $h$ of $H_{\Cal U}$ such that 
$x(g)\in U_h$. \endproclaim 
 
\demo{Proof} Define $x$ on $D_x\cap G_k$ inductively. First, $x$ is already 
defined on $D_x\cap G_0 = \{ e\}$: $x(e) = p$. Then, if $g\in D_x\cap G_k$, 
$x(g)$ is already defined, $h = h_{U,V}$ is a generator of $H_{\Cal U}$ and 
$x(g)$ lies in the domain of $h$, then join $x(g)$ to $h(x(g))$ by the 
unique leaf geodesic $\gamma : [0,1] \to W = W_{U,V}$, built from the fibres 
of the bundle $E$ a tubular neighbourhood $E_\gamma$ of a leaf $L_p$ around 
$\gamma ([0, 1])$  and lift $\gamma$ to the unique curve $\tilde\gamma : 
[0, 1]\to E_\gamma$ such that $\pi_\gamma\circ\tilde\gamma = \gamma$, 
$\tilde\gamma (0) = p$ and $\tilde\gamma (t)\in N$ for all $t\in [0,1]$, 
$\pi_\gamma$ being the canonical projection in the 
tubular neighbourhood. Finally, put $x(h, g) = \tilde\gamma (1)$
(see Figure 1 with $g = e$).

 
Since the lengths of all such geodesics $\gamma$ are uniformly bounded, this  
construction can be performed for $\alpha$ small enough, say 
$\alpha\le\alpha_0$. Also,  
the same reason together with uniform bounds for the geometry (i. e., sectional
curvature, injectivity radius etc.) of $M$ implies that $d_M(x(h, g), hx(g))
\le C\alpha$ and that $l(\tilde\gamma )\le C'l(\gamma )$ for some universal 
constants $C$ and $C'$.  
Therefore, this construction leads to an $C\alpha$-pseudo-orbit of $H_{\Cal U}$
which satisfies all the required conditions.  
 
The uniqueness of such an $x$ follows from the fact that at each step of 
induction the construction is in fact unique if $\alpha$ is as small that the 
projections $\pi_\gamma |N$ are one-to-one for all minimal geodesics joining 
$q$ to $h(q)$ for any $q\in M$ and any holonomy map $h$ generating 
$H_{\Cal U}$. \qed 
\enddemo


\demo{Remark} Clearly, any orbit of the holonomy pseudogroup $\Cal H_{\Cal U}$
is contained in a leaf of $\Cal F$. It seems that embedding a holonomy
$\alpha$-pseudo-orbit $x$ into a pseudoleaf is not so obvious. The difficulty
comes from the fact that, even if $\alpha$ is very small, the points $x(h_1)$
and $x(h_2)$, $h_1, h_2\in (\Cal H_{\Cal U})_1$, may belong to a disc
orthogonal to $\Cal F$ (see Figure 2, where $q = h_1(p) = h_2(p)$ and
$p_i = x(h_i)$) and therefore any submanifold $N$ containing
these points in a reasonable way has to make an angle close to $\frac{\pi}{2}$
with some leaves of $\Cal F$. However, in the case of foliated bundles,
one can produce [{\bf In}] immersed pseudoleaves corresponding to pseudo-orbits
of the global holonomy group by suitable use of centres of mass [{\bf BK}] 
of the points of pseudo-orbits.
\enddemo
 
\subhead 4. Geometric entropy\endsubhead 
Let again $(M, \Cal F, \langle\cdot , \cdot\rangle )$ be a compact foliated 
Riemannian manifold, $\Cal U$ - a nice cover of $(M, \Cal F)$ and $T = 
T_{\Cal U}$ - a complete transversal. Recall that two points $p, q\in T$ are 
$(R, \epsilon )$-{\it separated} w. r. t. $\Cal F$ if either 
$d_M(p, q)\ge\epsilon$ or there exists a leaf curve $\gamma :[0, 1]\to M$ of 
length $l(\gamma )\le R$ and such that $\gamma (0) = p$ (resp., $\gamma (0) 
= q$), the originated at an $\Cal U$-plaque through $q$ (resp., through $p$) 
orthogonal projection $\tilde\gamma$ of $\gamma$ onto the leaf through $q$ 
(resp., through $p$) exists and satisfies the condition 
$$d_M(\gamma (1), \tilde\gamma (1))\ge\epsilon ,\tag15$$ 
$d_M$ being, as before, the Riemannian distance function on $M$. As usually, a
subset $A$ of $T$ is $(R, \epsilon )$-separated when any two its points have 
this property. Since $M$ is compact, separated sets are always finite and we 
may define $N(R, \epsilon , \Cal F)$ as the maximal cardinality of an $(R, 
\epsilon)$-separated (w. r. t. $\Cal F$) set. The {\it geometric entropy} 
$h(\Cal F)$ of $\Cal F$ is defined by (see [{\bf GLW}] and [{\bf LW2}]) 
$$h(\Cal F) = \lim_{\epsilon\to 0} N(\epsilon , \Cal F), \tag16$$ 
where 
$$N(\epsilon , \Cal F) = \limsup_{R\to\infty}\frac{1}{R}\log N(R, \epsilon , 
\Cal F).\tag17$$ 
Obviously, $h(\Cal F)$ depends on the Riemannian metric $\langle\cdot , 
\cdot\rangle$ on $M$ or, more precisely, on the Riemannian metric 
$\langle\cdot , \cdot\rangle |T\Cal F\otimes T\Cal F$ along the leaves: 
If two Riemannian metrics coincide on  
$T\Cal F$, then the corresponding entropies are equal. 
 
It is well known that the geometric entropy $h(\Cal F)$ is strongly related 
to the entropies of the holonomy pseudogroups: 
$$h(\Cal F) = \sup_{\Cal U}\frac{1}{\Delta (\Cal U)}h(H_{\Cal U}),\tag18$$ 
where $\Cal U$ ranges over the family of all nice covers of $(M, \Cal F)$ 
and $\Delta (\Cal U)$ is the smallest upper bound for the diameters of the 
plaques of $\Cal U$ ([{\bf GLW}], Thm. 3.4).  
 
Note that, since any point $p$ determines the leaf $L_p$ uniquely, we can  
consider $N(R, \epsilon , \Cal F)$ as the number of separated leaves $L$ 
with base points $p\in T\cap L$. This approach can be adapted to the case 
of pseudoleaves even if there are infinitely many pseudoleaves passing 
through a given point. 
 
Fix $\alpha_0 > 0$. Again since $M$ is compact and has bounded geometry, there  
exist positive numbers $\epsilon_0$ and $\epsilon_1$ such that for any points 
$p_1, p_2\in M$ for which $d_M(p_1, p_2) < \epsilon_0$ and any 
$\alpha$-pseudoleaves $N_1$ and $N_2$ with $\alpha\le\alpha_0$ and $p_i\in 
N_i$ there exists an unique vector $v\in T^\perp_{p_1}N_1$ such that $p'_1 = 
\exp_M(v)\in N_2$ and $d_{N_2}(p'_1, p_2) < \epsilon_1$. The point $p'_1$ can 
be considered as the orthogonal projection of $p_1$ onto $N_2$. 
 
So, let us say that $\alpha$-pseudoleaves $(N_i, p_i)$ 
($\alpha\le\alpha_0$, $i = 1, 2$) with base  
points $p_i\in N_i\cap T$ are $(R, \epsilon )$-separated if either $d_M(p_1,  
p_2)\ge\epsilon_0$ or there exists a curve $\gamma : [0, 1]\to N_1$ (resp.,, 
$N_2$) originated at $p_1$ (resp., at $p_2$) of the length $l(\gamma )\le R$ 
and such that its orthogonal projection $\tilde\gamma$ onto $N_2$ (resp., 
onto $N_1$) originated in the ball $B_{N_2}(p_2, \epsilon_1)$ (resp., in 
$B_{N_1}(p_1, \epsilon_1)$) exists and satisfies inequality (15). Analogously 
to (16) and (17), let us put 
$$h_{\ps}(\Cal F) = \lim_{\epsilon\to 0^+}\inf_{\alpha} N_\alpha (\epsilon ,  
\Cal F),\tag18$$  
where 
$$N_\alpha (\epsilon , \Cal F) = \limsup_{R\to\infty}\frac{1}{R}\log 
N_{\alpha (R)}(R, \epsilon , \Cal F),\tag19$$ 
$\alpha :\Bbb R_+\to\Bbb R_+$ is any function satisfying $\alpha (R)\to 0$ 
while $R\to\infty$, and $N_{\alpha (R)}(R, \epsilon , \Cal F)$ is the maximal 
number of pairwise $(R, \epsilon )$-separated $\alpha (R)$-pseudoleaves 
equipped with base points. 
 
Since any leaf becomes an $\alpha$-pseudoleaf for any $\alpha\ge 0$, it is 
obvious that our "geometric pseudo-entropy" $h_{\ps}(\Cal F)$ is not less 
than the "real" geometric entropy  
$h(\Cal F)$. Below, we shall prove that, similarly as in the case of a 
pseudogroup, these two entropies are equal. 
 
\proclaim{Theorem 2} For any compact foliated Riemannian manifold $(M, \Cal F,
\langle\cdot , \cdot\rangle )$ the equality 
$$h_{\ps}(\Cal F) = h(\Cal F)\tag20$$ 
holds. 
\endproclaim 
 
\demo{Proof} As was observed, we have only to prove the inequality "$\le$" 
in (20). 
 
To this end cover $M$ by "coins" as in the proof of Theorem 3.4 of [{\bf GLW}].  
"Coins" are distinguished charts $U_p$, $p\in M$, built along the transversals 
$T_p = \exp_M B^\perp (0_p, \rho )$, $B^\perp (0_p, \rho )$ being the centred 
at the origin ball of radius $\rho$ in $T_p^\perp\Cal F$, the orthogonal 
complement of $T_p\Cal F$ in $T_pM$. More precisely, 
$$U_p = \cup_{q\in T_p} B^{\Cal F}(q, \Delta ),$$ 
where $B^{\Cal F}(q, \Delta )$ is a centred at $q$ ball of radius $\Delta$ 
in $L_q$, the leaf of $\Cal F$ at $q$. Here, $\rho$ and $\Delta$ are 
sufficiently small such that the corresponding exponential maps are 
diffeomorphic on all the balls of radii $2\rho$ and $2\Delta$.
 
Choose a finite subset $Q$ of $M$, $Q = \{ q_1,\dots , q_N\}$, such that 
the corresponding charts $U_{q_i}$, $i = 1, \dots , N$, form a nice  
covering $\Cal U_\Delta$. Obviously, $\Delta (\Cal U_\Delta ) = 2\Delta$. 
Also, fix a small number $\eta > 0$ and choose $Q$  dense in $M$ so that 
any point of $M$ is at most $\eta/2$ apart in the leaf distance from 
a point of $T = T_{\Cal U_\Delta}$ and  
$\eta/2$ apart in $d_M$-distance from one of the points $q_j\in Q$. 
The existence of such a set $Q$ is shown in full detail in [{\bf GLW}].
 
Take $\alpha > 0$ small enough, $\epsilon\in (0, \epsilon_0)$,  
$(R, \epsilon)$-separated $\alpha$-pseudoleaves $(N_i, p_i)$, $i = 1, 2$, 
and the corresponding $C\alpha$-pseudo-orbits $x_i$ of  $H_{\Cal U_\Delta}$ 
(Lemma 6).  
 
If $d_M(p_1, p_2)\ge\epsilon_0$, then $d(x_1(e), x_2(e))\ge\epsilon_0
\ge\epsilon$ and the pseudo-orbits $x_1$ and $x_2$ are $(n, 
\epsilon )$-separated for all $n\ge 1$.  
Otherwise, there exists a curve $\gamma :[0, 1]\to N_1$, $\gamma (0) = p_1$,  
$l(\gamma )\le R$, separating $N_1$ from $N_2$, i.e. such that 
$d(\gamma (1), \tilde\gamma (1))\ge\epsilon$, where $\tilde\gamma$ is the 
suitable projection of $\gamma$ onto $N_2$. Shortening $\gamma$ if necessary, 
we may assume that  $d(\gamma (1), \tilde\gamma (1)) = \epsilon$ and 
$d(\gamma (t), \tilde\gamma  (t))\le\epsilon$ for all $t$.

Split $\gamma$ into $n = [1 + R/(2\Delta - 2\eta )]$ pieces $\gamma_j$ of 
length $l(\gamma_j)\le 2\Delta - 2\eta$. Set $p'_0 = p_1$ and $p'_{n+1} 
= \gamma (1)$, for each $j = 1,\dots , n$ denote by $p'_j$  the midpoint of 
$\gamma_j$ and find indices $i(j)\in\{ 1,\dots , N\}$  such that $d_M(p'_j, 
q_{i(j)})\le\eta /2$ for $j = 0,1,\dots , 
n+1$. Thus, if $\alpha$ is small enough, each piece $\gamma_j$ is contained 
in a chart $U_j =  U_{q_{i(j)}}\in\Cal U_\Delta$. Let  $g = (h_n,\dots , h_0)$,
where  $h_j = h_{U_jU_{j+1}}$ for $j = 0,1,\dots , n$. 


From the construction of the pseudo-orbits $x_i$ it follows that $d_{N_1}(p'_j, 
x_1(g_j))\le\eta$ when $g_j = (h_j, \dots , h_0)$, $j = 0, 1, \dots , n$ and $\alpha$ is, 
again, small enough (see Figure 3, where $r_j = x_1(g_j)$, $r_{j+1} =
x_1(g_{j+1})$, $s_j = x_2(g_j)$ and $s_{j+1} = x_2(g_{j+1})$).
It follows that $x_1(g_j)\in U_j$ for all $j$ if only
$\eta$ is small with respect to $\Delta$ and $\rho$. In particular, $g\in D_{x_1}$. 
Also, $g\in D_{x_2}$ if only $\epsilon$ is small enough. Moreover, since 
$d_{N_1}(x_1(g), \gamma (1))\le\eta$ and the similar inequality holds for 
$\tilde\gamma (1)$ and $x_2(g)$ on $N_2$, then $d_M(x_1(g), x_2(g))\ge 
C''\epsilon$, where $C'' > 0$ is an universal constant depending only on $\rho$, 
$\Delta$, $\eta$ and the geometry of $M$. Hence, the pseudo-orbits  
$x_1$ and $x_2$ of $H_{\Cal U_\Delta}$ are $(n+1, C''\epsilon )$-strongly 
separated and consequently
$$N_\alpha (R, \epsilon , \Cal F)\le N_{C\alpha} ([2 + R/(2\Delta - 2\eta)],  
C''\epsilon , H_{\Cal U_\Delta})$$ 
for $\alpha$ and $\epsilon$ small enough.  
 
Now, pass to suitable limits to arrive at 
$$h_{\ps}(\Cal F)\le\frac{1}{2\Delta - 2\eta}h_{\ps}(H_{\Cal U_\Delta}).$$ 
Finally, apply Theorem 1 and  relation (18) to obtain the inequality 
$$h_{\ps}(\Cal F)\le\frac{\Delta}{\Delta - \eta}h(\Cal F)$$ 
which yields the required result since $\eta > 0$ may be arbitrarily small. \qed 
\enddemo 

\subhead 5. Groups and foliated bundles\endsubhead 
If $G$ is a group of  global transformations of $X$, then one need not take care of
domains of maps, $D_x = G_\infty$ for any pseudo-orbit $x$,  $d_0 = d_1$ is defined  
simply by 
$$d_0(x, y) = d_1(x, y) = \sum_{k=0}^\infty\frac{1}{2^k}\max\{ d(x(g), y(g)); g\in  
G_k\}$$ 
and the maps $\sigma_g : Y_\alpha\to Y_\alpha$, $g\in G_\infty$ are also defined  
globally. In other words, if the group $G$ acts on $X$, then it acts on the spaces of  
pseudo-orbits as well. 
 
In particular, if $\Cal F$ is a foliation transverse to the fibres of a fibre bundle $\pi  
:M\to B$, then its global holonomy group $H$ acts on a fibre $F$ (over a fixed point  
$b_0\in B$). $H$ is generated by the maps $h_{[\gamma ]}$ with $[\gamma ]\in\pi_1  
(B, b_0)$ generating the fundamental group of $B$, defined by 
$$h_{[\gamma ]}(p) = \tilde\gamma (1),$$ 
where $\tilde\gamma :[0, 1]\to L_p$ is the originated at $p$ lift of $\gamma$ 
to the leaf $L_p$ through $p\in F$. 
 
Similarly, if $N$ is an $\alpha$-pseudoleaf and $\alpha$ is small enough, then $N$ is  
transverse to the fibres of our bundle and the corresponding  pseudo-orbit of $H$ can  
be obtained  by lifting to $N$ the loops generating $\pi_1(B, b_0)$.  
 
For such a foliated bundle, a Riemannian metric on $M$ which, when restricted to  
$T\Cal F$, coincides with the lift of a Riemannian metric on $B$, and a generating set  
$H_1\subset H$ one has the inequality 
$$\frac{1}{a}h(H, H_1)\le h(\Cal F)\le\frac{1}{b}h(H, H_1),\tag21$$ 
where $a$ (resp. $b$) is the maximum (resp., minimum) of the lengths of the  
homotopy classes $[\gamma ]\in\pi_1(B, b_0)$ with $h_{[\gamma ]}\in H_1$ ([{\bf  
GLW}], p. 122). In the same way, 
$$\frac{1}{a}h_{\ps}(H, H_1)\le h_{\ps}(\Cal F)\le\frac{1}{b}h_{\ps}(H,  
H_1),\tag22$$ 
Since $h_{\ps}(H, H_1) = h(H, H_1)$ by Theorem 1, inequalities (21) and (22) yield 
$$h_{\ps}(\Cal F)\le\frac{a}{b}h(\Cal F).\tag23$$ 
 
Fix $m\in\Bbb N$ and replace $H_1 = \{z_1, \dots z_s\}$ by a new generating set  
$$H'_1 = \{ z_1^{m_1}, z_1^{m_1+1},\dots , z_s^{m_s}, z_s^{m_s+1}\}$$ 
with the exponents $m_k = [m/l_k]$, $l_k$ being the length of the homotopy class 
$z_k$. For this generating set, the corresponding factor $a'/b'$ in (23) satisfies the 
inequalities
$$\frac{m}{m + 1 + b}\le\frac{a'}{b'}\le\frac{m + 1 + b}{m}$$
and is arbitrarily close to 1 when $m$ is sufficiently large.
 
This implies Theorem 2 for foliated bundles immediately. 

\subhead 6. An example\endsubhead
In [{\bf BS}] and [{\bf Mi}], the equality
$$h(G, G_1) = \lim_{\epsilon\to 0^+}\lim_{\alpha\to 0^+}\limsup_{n\to\infty}
\frac{1}{n}\log N_\alpha (n, \epsilon )\tag24$$
was proved for the group $G = \Bbb Z$ generated by $G_1 = \{ f\}$, a single
transformation of $X$. Reading the proof of (24) one can observe that some
of the argument is very specific: It depends strongly on the fact that
the growth of $\# G_n$, $n\in\Bbb N$, is linear. Below, we show an example,
where the equality (24) does not hold: The limit on the right hand side
is strongly bigger than the entropy.

So, let $X = S^1$, $G = F_2$ be the free group with two generators,  
$f_0$ and $ f_1$, $f_i$ being Morse - Smale diffeomorphisms of $S^1$, each of 
them with exactly two fixed points, a source $p_i$ and a sink $q_i$, such that
$\{p_0 , q_0\}\cap\{ p_1, q_1\} = \emptyset$ and the sets $U_i = \{ p\in S_1,
f'_i > 1 + \delta\ \text{on}\ [p - \alpha_0, p + \alpha_0]\}$, $i = 0, 1$,
cover $S^1$ for some $\alpha_0$ and $\delta > 0$. 
Also, let $G_1 = \{\id , f_0^{\pm 1}, f_1^{\pm 1}\}$.

Now, take any $n\in\Bbb N$, $\epsilon > 0$, $\alpha\in (0, \alpha_0)$ and an 
$(n, \epsilon )$-separated set $A\subset S^1$. For any $x\in A$ 
and any $g = (g_n, \dots , g_1)\in G_n$ with $g_i\in\{ f_0, f_1\}$ 
define an $\alpha$-pseudo-orbit $\tilde x_g$ of $G$ as follows: 
\roster
\item $\tilde x_g (h) = h(x)$ for all $h\in G_n$ and any $h\in G_m\setminus 
G_{m-1}$, $m > n$,
of the form $h = (g', g'')$ with $g'\in G_{m-n}$ and $g''\in G_n$, $g''\ne g$.
\item If $g(x)\in U_0$, then $h_1 = f_0$, otherwise $h_2 = f_1$, and, in both
cases, $\tilde x_g (h_1, g) = h_1(g(x)) + \alpha$.
\item If $h_1, \dots , h_j\in G_1$ and $\tilde x_g (h_j,\dots , h_1, g)$ are 
already defined, and \linebreak $\tilde x_g(h_j,\dots , h_1, g)\in U_0$, then
$h_{j+1} = f_0$, otherwise $h_{j+1} = f_1$, and, in both cases,
$$\tilde x_g (h_{j+1}, h_j, \dots , h_1, g) = h_{j+1}(\tilde
x_g (h_j,\dots , h_1, g)) + \alpha .$$
\item If $j\in\Bbb N$, $m\in\Bbb N$, $h = (h'_m,\dots, h'_1, h_j,\dots, h_1, g)$
with $h'_i\in G_1$, $h'_1\ne h_{j+1}$ and $h'_1\ne h_j^{-1}$, then
$\tilde x_g (h) = h'_m\circ\dots\circ h'_1(\tilde x_g(h_j, \dots , h_1, g))$.
\endroster
Roughly speaking, $\tilde x_g$ is obtained by a suitable modification of the 
orbit $G(x)$ along a single originated at $g$ branch of the Cayley graph of $G$.

A standard induction involving Mean Value Theorem shows that 
$$d(\tilde x_g(h_j, \dots h_1, g), h_j\circ\dots\circ h_1\circ g(x))
\ge\alpha\cdot (1 + \delta )^j\ge\epsilon\tag25$$
whenever $j > n(\epsilon , \alpha )$, $n(\epsilon , \alpha )$ being the 
largest integer which does not exceed the quotient
$$\log\frac{\epsilon}{\alpha}/\log (1 + \delta).$$
It follows that the set $A_\alpha$ of all such pseudo-orbits $x_g$ is 
$(n + n(\epsilon , \alpha) + 1, \epsilon)$-separated. In fact, if $x\ne x'\in A$,
then $d(\tilde x_g(h), \tilde x'_{g'}(h)) = d(h(x), h(x'))\ge\epsilon$ for some 
$h\in G_n$; if $x = x'$ but $g\ne g'$, then $\tilde x_{g'}(h_j, \dots , h_1, g) = 
h_j\circ\dots\circ h_1\circ g(x)$ and the inequality $d(\tilde x_{g'}(h), 
\tilde x_g(h))\ge\epsilon$ for some $h\in G_{n + n(\epsilon , \alpha ) + 1}$ 
follows from (25). Clearly,
$$\# A_\alpha = 2^n\cdot\# A.$$
This yields the inquality
$$N_\alpha (n + n(\epsilon , \alpha) + 1, \epsilon )\ge 2^n\cdot N(n, 
\epsilon , S^1)$$
for all $n$, $\epsilon > 0$ and $\alpha < \alpha_0$. Passing to suitable limits
we end up with
$$\lim_{\epsilon\to 0^+}\lim_{\alpha\to 0^+}\limsup_{n\to\infty}
\frac{1}{n}\log N_\alpha (n, \epsilon )\ge h(G, G_1) + \log 2.\tag26$$

Inequality (26) shows that, what was already observed in [{\bf GLW}], the 
entropy of a group (a pseudogroup, a foliation) depends strongly on the growth 
of the group (the orbits, the leaves). It would be very interesting to obtain 
a good upper estimate of the right hand side in (24) in terms of $h(G, G_1)$ and
the growth of $G$. Some indication in this direction is contained in the proof
of Proposition 6.13 in [{\bf GLW}]. 


\Refs 
\widestnumber\key{\bf GLW} 
\ref\key{\bf BS}\by M. Barge and R. Swanson\paper Pseudo-orbits and topological  
entropy\jour Proc. Amer. Math. Soc.\vol 109\yr 1990\pages 559--566\endref 

\ref\key{\bf BK}\by P. Buser and H. Karcher\paper Gromov's almost flat 
manifolds\jour Ast\'erisque\vol 81\yr 1981\pages 1--148\endref
  
\ref\key{\bf Eg}\by S. Egashira\paper Expansion growth of foliations\jour Ann. Fac.  
Sci. Univ. Toulouse\vol 2\yr 1993\pages 15--52\endref 
 
\ref\key{\bf GLW}\by E. Ghys, R. Langevin and P. Walczak\paper Entropie  
g\'eom\'etrique des feuilletages\jour Acta Math.\vol 160\yr 1988\pages 105--142\endref 
 
\ref\key{\bf HH}\by G. Hector and u. Hirsch\book Introduction to the Geometry of  
Foliations, Part A and B\publ Vieweg and Sohn\publaddr Braunschweig\yr 1986 and  
1987\endref 
 
\ref\key{\bf Hu1}\by S. Hurder\paper Ergodic theory of foliations and a theorem of  
Sacksteder \inbook Dynamical Systems: Proceedings, University of Maryland 1986 -  
87. Lecture Notes in Math., vol 1342\publ Springer Verlag\publaddr Berlin and New  
York\yr 1988 \pages 291--328\endref 
 
\ref\key{\bf Hu2}\by S. Hurder\paper Exceptional minimal sets of C$^{1+\alpha}$  
group actions on the circle\jour Ergodic Th. \& Dynam. Sys.\vol 11
\yr 1991\pages 455--467\endref 
 
\ref\key{\bf Hu3}\by S. Hurder\paper Coarse geometry of foliations\inbook Geometric  
Study of Foliations, Proceedings, Tokyo 1993\publ World Sci.\publaddr  
Singapore\pages 35--96 \endref 
 
\ref\key{\bf Hur}\by M. Hurley\paper On topological entropy of maps\jour Ergodic 
Th. \& Dynam. Sys.\vol 15\yr 1995\pages 557--568\endref 
 
\ref\key{\bf In}\by T. Inaba\paper Expansivity, pseudoleaf tracing property and  
topological stability of foliations\paperinfo Preprint 1995\endref 
 
\ref\key{\bf IT}\by T. Inaba and N. Tsuchiya\paper Expansive foliations\jour 
Hokkaido Math. J.\vol 21\yr 1992\pages 39--49\endref 
 
\ref\key{\bf LW1}\by R. Langevin and P. Walczak\paper Entropie d'une  
dynamique\jour C. R. Acad. Sci. Paris, S\'er. I\vol 312\yr 1991\pages 141--144\endref 
 
\ref\key{\bf LW2}\by R. Langevin and P. Walczak\paper Entropy, transverse entropy  
and partitions of unity\jour Ergodic Th. \& Dynam. Sys.\vol 14\yr 1994\pages 551-- 
563\endref 
 
\ref\key{\bf Mi}\by M. Misiurewicz\paper Remark on the definition of topological  
entropy\inbook Dynamical Systems and Partial Differential Equations (Caracas  
1984)\publaddr Caracas\yr 1986\pages 65--67\endref 
\endRefs
\enddocument